\documentclass[11pt,letterpaper]{article}
\RequirePackage[backend=biber,bibstyle=authoryear,citestyle=authoryear,
natbib=true,sorting=nyt,sortcites=true,maxnames=10,minnames=10,
hyperref=true,abbreviate=false,date=long,isbn=true,eprint=true,
uniquename=init,giveninits=true]{biblatex}
\usepackage{hyphenat,graphicx}
\usepackage[colorlinks,linktocpage,breaklinks]{hyperref}
\makeatletter
\def\@filecolor{blue}
\def\@linkcolor{blue}
\def\@citecolor{blue}
\def\@urlcolor{blue}
\let\@old@citep\citep
\let\@old@citet\citet
\let\@old@citeauthor\citeauthor
\def\citep{\@old@citep*}
\def\citet{\@old@citet*}
\def\citeauthor{\@old@citeauthor*}
\let\cite\citep
\makeatother
\usepackage[top=1.5in,bottom=1.2in,left=1.35in,right=1.35in,letterpaper]%
  {geometry}
\makeatletter
\def\ps@headings{%
    \let\@mkboth\@gobbletwo
    \def\@oddhead{\hss\scshape\shorttitle\hss\reset@font\rmfamily\thepage}
    \def\@evenhead{\reset@font\rmfamily\thepage\hss\scshape\shortauthors\hss}
    \let\@oddfoot\@empty\let\@evenfoot\@empty}
\@twosidetrue
\makeatother
\usepackage[neveradjust]{paralist}
\usepackage{theorem}
\newtheorem{theorem}{Theorem}[section]
\makeatletter
\@namedef{procnonum}{\trivlist \@topsep \theorempreskipamount
  \@topsepadd \theorempostskipamount
  \@ifnextchar[{\@ADL@xprocnonumstar}{\@ADL@yprocnonumstar}}
\@namedef{endprocnonum}{\endtrivlist}
\def\@ADL@xprocnonumstar[#1]{\item[\hskip \labelsep{\theorem@headerfont #1}]
  \normalfont\rmfamily\itshape}
\def\@ADL@yprocnonumstar{\item[] \normalfont\rmfamily\itshape}
\makeatother
\makeatletter
\newcommand\pushright{\protect\@ADL@pushright}
\newcommand\@ADL@pushright[1]{{\ifvmode\null\hfill{#1}\par\else\ifmmode%
  \@ADLmaths@pushright{\hbox{#1}}\else\ifinner\@ADLhbox@pushright{#1}%
  \else\@ADLparag@pushright{#1}\fi\fi\fi}}
\newcommand\@ADLmaths@pushright[1]{{\ifinner\@ADLhbox@pushright{#1}\else%
  \tag*{$#1$}\fi}}
\newcommand\@ADLparag@pushright[1]{{\parfillskip=0pt\widowpenalty=10000%
  \displaywidowpenalty=10000\finalhyphendemerits=0\@ADLhbox@pushright#1\par}}
\newcommand\@ADLhbox@pushright{\unskip\nobreak\hfil\penalty50\hskip.2em%
  \null\hfill\hfill}
\newenvironment{proof}{\trivlist\item[\hskip\labelsep\textit{Proof:}\/]%
  \@ADLsave@set@qed\xspace\normalfont\rmfamily}
  {\qed\@ADLrestore@qed\endtrivlist}
\newenvironment{subproof}{\trivlist\item[\hskip\labelsep\textit{Proof:}\/]%
  \@ADLsave@set@subqed\normalfont\rmfamily}
  {\subqed\@ADLrestore@subqed\endtrivlist}
\newif\if@ADL@qed\@ADL@qedfalse
\newcommand\qed{\protect\@ADL@qed{$\blacksquare$}}
\newcommand\@ADL@qed[1]{\if@ADL@qed\global\@ADL@qedfalse%
  \pushright{#1}\else\ifhmode\ifinner\else\par\fi\fi\fi}
\newcommand\@ADLrestore@qed{\global\let\if@ADL@qed\@ADLsaved@ifqed}
\newcommand\@ADLsave@set@qed{\let\@ADLsaved@ifqed
  \if@ADL@qed\global\@ADL@qedtrue}

\newif\if@ADL@subqed\@ADL@subqedfalse
\newcommand\subqed{\protect\@ADL@subqed{$\blacktriangledown$}}
\newcommand\@ADL@subqed[1]{\if@ADL@subqed\global\@ADL@subqedfalse%
  \pushright{#1}\else\ifhmode\ifinner\else\par\fi\fi\fi}
\newcommand\@ADLrestore@subqed{\global\let\if@ADL@subqed\@ADLsaved@ifsubqed}
\newcommand\@ADLsave@set@subqed{\let\@ADLsaved@ifsubqed
  \if@ADL@subqed\global\@ADL@subqedtrue}
  
\makeatother  
\usepackage[reqno,centertags]{amsmath}
\usepackage{amssymb,xspace}
\makeatletter
\def\interval{\@ifnextchar({\@ADL@openleftint}{\@ADL@closedleftint}}
\def\@ADL@openleftint(#1,#2{(#1,#2%
  \@ifnextchar){\@ADL@openrightint}{\@ADL@closedrightint}}
\def\@ADL@closedleftint[#1,#2{[#1,#2%
  \@ifnextchar){\@ADL@openrightint}{\@ADL@closedrightint}}
\def\@ADL@openrightint){)}
\def\@ADL@closedrightint]{]}
\makeatother
\newenvironment{keywords}{\quote\small\textbf{Keywords.}}{\endquote}
\newenvironment{AMS}{\quote\small\textbf{AMS Subject Classifications (2010).}}
   {\endquote}
\newcommand\defn[1]{{\normalfont\bfseries\emph{\mathversion{bold}#1}}}
\usepackage{booktabs,diagbox,pict2e}
\newcommand\integer{\mathbb{Z}}
\newcommand\integerp{\integer_{>0}}
\newcommand\real{\mathbb{R}}
\newcommand\realp{\real_{>0}}
\newcommand\realnn{\real_{\ge0}}

\newcommand\vect[1]{\boldsymbol{#1}}
\newcommand\mat[1]{\boldsymbol{#1}}
\newcommand\transpose[1]{#1^T}
\newcommand\sign{\operatorname{sign}}
\newcommand\subsets{\boldsymbol{2}}
\newcommand\card{\operatorname{card}}
\newcommand\ie{i.e.,}
\newcommand\eg{e.g.,}

\newcommand\scirc{\raise1pt\hbox{$\,\scriptstyle\circ\,$}}
\newcommand\map[3]{#1\colon#2\rightarrow#3}
\newcommand\mapdef[5]{\begin{aligned}
  #1\colon&\begin{aligned}[t]#2\end{aligned}\rightarrow
  \begin{aligned}[t]#3\end{aligned}\\&\begin{aligned}[t]#4\end{aligned}
  \mapsto\begin{aligned}[t]#5\end{aligned}\end{aligned}}
  \newcommand\lin[2]{\textup{L}(#1;#2)}
\makeatletter
\newcommand\linder{\@ifnextchar[{\@ADL@rlinder}{\@ADL@linder}}
\def\@ADL@rlinder[#1]#2{\boldsymbol{D}^{#1}#2}
\newcommand\@ADL@linder[1]{\boldsymbol{D}#1}
\newcommand\plinder{\@ifnextchar[{\@ADL@rplinder}{\@ADL@plinder}}
\def\@ADL@rplinder[#1]#2#3{\boldsymbol{D}^{#1}_{#2}#3}
\newcommand\@ADL@plinder[2]{\boldsymbol{D}_{#1}#2}
\makeatother
\newcommand\slnorm{\lvert}
\newcommand\srnorm{\rvert}
\newcommand\snorm[1]{\slnorm #1\srnorm}
\newcommand\asnorm[1]{\left\slnorm #1\right\srnorm}
\newcommand\dlnorm{\lVert}
\newcommand\drnorm{\rVert}
\newcommand\dnorm[1]{\dlnorm #1\drnorm}
\newcommand\adnorm[1]{\left\dlnorm #1\right\drnorm}
\newcommand\tnorm[1]{|||#1|||}
\newcommand\setdef[2]{\{#1\;|\enspace#2\}}
\newcommand\asetdef[2]{\left\{#1\immediate\vphantom{#2}\;\right|
  \left.\immediate\vphantom{#1}\enspace#2\right\}}
\newcommand\inprod[2]{\langle#1,#2\rangle}
\title{A top nine list: Most popular induced matrix norms}
\author{Andrew D.\ Lewis\thanks{Professor, Department of Mathematics and
Statistics, Queen's University, Kingston, ON K7L 3N6, Canada,
email:~\texttt{andrew.lewis@queensu.ca}}}
\date{2010/03/20}
\newcommand\shorttitle{A top nine list: Most popular induced matrix norms}
\newcommand\shortauthors{A.\ D.\ Lewis}
\begin{document}
\maketitle

\begin{abstract}
Explicit formulae are given for the nine possible induced matrix norms
corresponding to the $1$-, $2$-, and $\infty$-norms for Euclidean space.  The
complexity of computing these norms is investigated.
\end{abstract}


\section{Introduction}

Arguably the most commonly used norms for real Euclidean space $\real^n$ are
the norms $\dnorm{\cdot}_1$\@, $\dnorm{\cdot}_2$\@, and
$\dnorm{\cdot}_\infty$ defined by
\begin{equation*}
\dnorm{\vect{x}}_1=\sum_{j=1}^n\snorm{x_j},\quad
\dnorm{\vect{x}}_2=\left(\sum_{j=1}^n\snorm{x_j}^2\right)^{1/2},\quad
\dnorm{\vect{x}}_\infty=\max\{\snorm{x_1},\dots,\snorm{x_n}\},
\end{equation*}
respectively, for $\vect{x}=(x_1,\dots,x_n)\in\real^n$\@.  Let
$\lin{\real^n}{\real^m}$ be the set of linear maps from $\real^n$ to
$\real^m$\@, which we identify with the set of $m\times n$ matrices in the
usual way.  If $\mat{A}\in\lin{\real^n}{\real^m}$ and if
$p,q\in\{1,2,\infty\}$ then the norm of $\mat{A}$ induced by the $p$-norm on
$\real^n$ and the $q$-norm on $\real^m$ is
\begin{equation*}
\dnorm{\mat{A}}_{p,q}=\sup\setdef{\dnorm{\mat{A}(\vect{x})}_q}
{\dnorm{\vect{x}}_p=1}.
\end{equation*}
This is well-known to define a norm on $\lin{\real^n}{\real^m}$\@.  There are
other equivalent characterisations of the induced norm, but the one given
above is the only one we will need.  We refer to~\cite{RAH/CRJ:13} for a
general discussion of induced matrix norms.

For certain combinations of $(p,q)$\@, explicit expressions for
$\dnorm{\cdot}_{p,q}$ are known.  For example, in~\cite{RAH/CRJ:13}
expressions are given in the cases $(1,1)$ (in \S5.6.4), $(2,2)$ (\S5.6.6),
and $(\infty,\infty)$ (\S5.6.5).  In \cite{JR:00} the case $(\infty,1)$ is
studied, and its computation is shown to be NP-hard.  The case $(2,1)$ is
given by \citet{KD/BAP:09}\@, although the details of the degenerate case
given there are a little sketchy.  \citeauthor{KD/BAP:09} also list all of
the other combinations except $(2,\infty)$\@, for which no expression seems
to be available, and which we give here, apparently for the first time.  The
formula given by \citeauthor{KD/BAP:09} for $(\infty,2)$ is presented without
reference or proof, and is incorrect, probably a typographical error.

Here we present the correct formulae for all nine of the induced norms.
Although most of these formulae are known in the literature, we give proofs
in all nine cases so that, for the first time, all proofs for all cases are
given in one place.  We also analyse the computational complexity of
computing these various norms.

Here is the notation we use.  By $\{\vect{e}_1,\dots,\vect{e}_n\}$ we denote
the standard basis for $\real^n$\@.  For a matrix
$\mat{A}\in\lin{\real^n}{\real^m}$\@, $\vect{r}(\mat{A},a)\in\real^n$ denotes
the $a$th row and $\vect{c}(\mat{A},j)\in\real^m$ denotes the $j$th column.
The components of $\mat{A}$ are denoted by $A_{aj}$\@, $a\in\{1,\dots,m\}$\@,
$j\in\{1,\dots,n\}$\@.  The transpose of $\mat{A}$ is denoted by
$\transpose{\mat{A}}$\@.  The Euclidean inner product is denoted by
$\inprod{\cdot}{\cdot}$\@.  For a differentiable map
$\map{\vect{f}}{\real^n}{\real^m}$\@,
$\linder{\vect{f}}(\vect{x})\in\lin{\real^n}{\real^m}$ denotes the derivative
of $\vect{f}$ at $\vect{x}$\@.  For a set $X$\@, $\subsets^X$ denotes the
power set of $X$\@.

\section{Formulae for induced norms}

\begin{theorem}\label{the:matnorms}
Let\/ $p,q\in\{1,2,\infty\}$ and let\/ $\mat{A}\in\lin{\real^n}{\real^m}$\@.
The induced norm\/ $\dnorm{\cdot}_{p,q}$ satisfies the following formulae:
\begin{compactenum}[(i)]
\item \label{pl:matnorm11} $\dnorm{\mat{A}}_{1,1}=
\max\setdef{\dnorm{\vect{c}(\mat{A},j)}_1}{j\in\{1,\dots,n\}}$\@;
\item \label{pl:matnorm12} $\dnorm{\mat{A}}_{1,2}=
\max\setdef{\dnorm{\vect{c}(\mat{A},j)}_2}{j\in\{1,\dots,n\}}$\@;
\item \label{pl:matnorm1infty} $\displaystyle\begin{aligned}[t]
\dnorm{\mat{A}}_{1,\infty}=&\;
\max\setdef{\snorm{A_{aj}}}{a\in\{1,\dots,m\},\ j\in\{1,\dots,n\}}\\
=&\;\max\setdef{\dnorm{\vect{c}(\mat{A},j)}_\infty}{j\in\{1,\dots,n\}}\\
=&\;\max\setdef{\dnorm{\vect{r}(\mat{A},a)}_\infty}{a\in\{1,\dots,m\}}
\end{aligned}$\@;
\item \label{pl:matnorm21} $\dnorm{\mat{A}}_{2,1}=
\max\setdef{\dnorm{\transpose{\mat{A}}(\vect{u})}_2}{\vect{u}\in\{-1,1\}^m}$\@;
\item \label{pl:matnorm22} $\dnorm{\mat{A}}_{2,2}=
\max\setdef{\sqrt{\lambda}} {\lambda\ \textrm{is an
eigenvalue for}\ \transpose{\mat{A}}\mat{A}}$\@;
\item \label{pl:matnorm2infty} $\dnorm{\mat{A}}_{2,\infty}=
\max\setdef{\dnorm{\vect{r}(\mat{A},a)}_2}{a\in\{1,\dots,m\}}$\@;
\item \label{pl:matnorminfty1} $\dnorm{\mat{A}}_{\infty,1}=
\max\setdef{\dnorm{\mat{A}(\vect{u})}_1}{\vect{u}\in\{-1,1\}^n}$\@;
\item \label{pl:matnorminfty2} $\dnorm{\mat{A}}_{\infty,2}=
\max\setdef{\dnorm{\mat{A}(\vect{u})}_2}{\vect{u}\in\{-1,1\}^n}$\@;
\item \label{pl:matnorminftyinfty} $\dnorm{\mat{A}}_{\infty,\infty}=
\max\setdef{\dnorm{\vect{r}(\mat{A},a)}_1}{a\in\{1,\dots,m\}}$\@.
\end{compactenum}
\begin{proof}
\eqref{pl:matnorm11} We compute
\begin{align*}
\dnorm{\mat{A}}_{1,1}=&\;\sup\setdef{\dnorm{\mat{A}(\vect{x})}_1}
{\dnorm{\vect{x}}_1=1}\\
=&\;\sup\asetdef{\sum_{a=1}^m
\snorm{\inprod{\vect{r}(\mat{A}(\vect{x}))}{\vect{x}}}}
{\dnorm{\vect{x}}_1=1}\\
\le&\;\sup\asetdef{\sum_{a=1}^m\sum_{j=1}^n\snorm{A_{aj}}\snorm{x_j}}
{\dnorm{\vect{x}}_1=1}\\
=&\;\sup\asetdef{\sum_{j=1}^n\snorm{x_j}
\left(\sum_{a=1}^m\snorm{A_{aj}}\right)}{\dnorm{\vect{x}}_1=1}\\
\le&\;\max\asetdef{\sum_{a=1}^m\snorm{A_{aj}}}{j\in\{1,\dots,n\}}\\
=&\;\max\setdef{\dnorm{\vect{c}(\mat{A},j)}_1}{j\in\{1,\dots,n\}}.
\end{align*}
To establish the opposite inequality, suppose that $k\in\{1,\dots,n\}$ is
such that
\begin{equation*}
\dnorm{\vect{c}(\mat{A},k)}_1=
\max\setdef{\dnorm{\vect{c}(\mat{A},j)}_1}{j\in\{1,\dots,n\}}.
\end{equation*}
Then,
\begin{equation*}
\dnorm{\mat{A}(\vect{e}_k)}_1=\sum_{a=1}^m
\asnorm{\left(\sum_{j=1}^nA_{aj}\vect{e}_{k,j}\right)}=
\sum_{a=1}^m\snorm{A_{ak}}=\dnorm{\vect{c}(\mat{A},k)}_1.
\end{equation*}
Thus
\begin{equation*}
\dnorm{\mat{A}}_{1,1}\ge
\max\setdef{\dnorm{\vect{c}(\mat{A},j)}_1}{j\in\{1,\dots,n\}},
\end{equation*}
since $\dnorm{\vect{e}_k}_1=1$\@.

\eqref{pl:matnorm12} We compute
\begin{align*}
\dnorm{\mat{A}}_{1,2}=&\;\sup\setdef{\dnorm{\mat{A}(\vect{x})}_2}
{\dnorm{\vect{x}}_1=1}\\
=&\;\sup\asetdef{\left(\sum_{a=1}^m
\inprod{\vect{r}(\mat{A},a)}{\vect{x}}^2\right)^{1/2}}{\dnorm{\vect{x}}_1=1}\\
\le&\;\sup\asetdef{\left(\sum_{a=1}^m\left(\sum_{j=1}^n
\snorm{A_{aj}x_j}\right)^2\right)^{1/2}}{\dnorm{\vect{x}}_1=1}\\
\le&\;\sup\asetdef{\left(\sum_{a=1}^m
(\max\setdef{\snorm{A_{aj}}}{j\in\{1,\dots,n\}})^2
\left(\sum_{j=1}^n\snorm{x_j}\right)^2\right)^{1/2}}{\dnorm{\vect{x}}_1=1}\\
=&\;\left(\sum_{a=1}^m(\max\setdef{\snorm{A_{aj}}}
{j\in\{1,\dots,n\}})^2\right)^{1/2}\\
=&\;\left(\max\asetdef{\sum_{a=1}^mA_{aj}^2}
{j\in\{1,\dots,n\}}\right)^{1/2}=
\max\setdef{\dnorm{\vect{c}(\mat{A},j)}_2}{j\in\{1,\dots,n\}},
\end{align*}
using the fact that
\begin{equation*}
\sup\setdef{\dnorm{\vect{x}}_2}{\dnorm{\vect{x}}_1=1}=1.
\end{equation*}
To establish the other inequality, note that if we take $k\in\{1,\dots,n\}$
such that
\begin{equation*}
\dnorm{\vect{c}(\mat{A},k)}_2=\max\setdef{\dnorm{\vect{c}(\mat{A},j)}_2}
{j\in\{1,\dots,n\}},
\end{equation*}
then we have
\begin{equation*}
\dnorm{\mat{A}(\vect{e}_k)}_2=\left(\sum_{a=1}^m\left(\sum_{j=1}^n
A_{aj}\vect{e}_{k,j}\right)^2\right)^{1/2}=
\left(\sum_{a=1}^mA_{ak}^2\right)^{1/2}=\dnorm{\vect{c}(\mat{A},k)}_2.
\end{equation*}
Thus
\begin{equation*}
\dnorm{\mat{A}}_{1,2}\ge\max\setdef{\dnorm{\vect{c}(\mat{A},j)}_2}
{j\in\{1,\dots,n\}},
\end{equation*}
since $\dnorm{\vect{e}_k}_1=1$\@.

\eqref{pl:matnorm1infty} Here we compute
\begin{align*}
\dnorm{\mat{A}}_{1,\infty}=&\;\sup\setdef{\dnorm{\mat{A}(\vect{x})}_\infty}
{\dnorm{\vect{x}}_1=1}\\
=&\;\sup\asetdef{\max\asetdef{\asnorm{\sum_{j=1}^nA_{aj}x_j}}
{a\in\{1,\dots,m\}}}{\dnorm{\vect{x}}_1=1}\\
\le&\;\sup\asetdef{\max\asetdef{\snorm{A_{aj}}}{j\in\{1,\dots,n\},\
a\in\{1,\dots,m\}}\left(\sum_{j=1}^n\snorm{x_j}\right)}{\dnorm{\vect{x}}_1=1}\\
=&\;\max\setdef{\snorm{A_{aj}}}{j\in\{1,\dots,n\},\
a\in\{1,\dots,m\}}.
\end{align*}
For the converse inequality, let $k\in\{1,\dots,n\}$ be such that
\begin{equation*}
\max\setdef{\snorm{A_{ak}}}{a\in\{1,\dots,m\}}=
\max\setdef{\snorm{A_{aj}}}{j\in\{1,\dots,n\},\ a\in\{1,\dots,m\}}.
\end{equation*}
Then
\begin{align*}
\dnorm{\mat{A}(\vect{e}_k)}_\infty=&\;
\max\asetdef{\asnorm{\sum_{j=1}^nA_{aj}\vect{e}_{k,j}}}
{a\in\{1,\dots,m\}}\\
=&\;\max\setdef{\snorm{A_{ak}}}{a\in\{1,\dots,m\}}.
\end{align*}
Thus
\begin{equation*}
\dnorm{\mat{A}}_{1,\infty}\ge
\max\setdef{\snorm{A_{aj}}}{j\in\{1,\dots,n\},\ a\in\{1,\dots,m\}},
\end{equation*}
since $\dnorm{\vect{e}_k}_1=1$\@.

\eqref{pl:matnorm21} In this case we maximise the function
$\vect{x}\mapsto\dnorm{\mat{A}(\vect{x})}_1$ subject to the constraint that
$\dnorm{\vect{x}}_2^2=1$\@.  We shall do this using the Lagrange Multiplier
Theorem~\cite[\eg][\S{}II.5]{CHE:73}\@, defining
\begin{equation*}
f(\vect{x})=\dnorm{\mat{A}(\vect{x})}_1,\quad
g(\vect{x})=\dnorm{\vect{x}}_2^2-1.
\end{equation*}
Let us first assume that none of the rows of $\mat{A}$ are zero.  We must
exercise some care because $f$ is not differentiable on $\real^n$\@.
However, $f$ is differentiable at points off the set
\begin{equation*}
B_{\mat{A}}=\setdef{\vect{x}\in\real^n}{\textrm{there exists}\
a\in\{1,\dots,m\}\ \textrm{such that}\
\inprod{\vect{r}(\mat{A},a)}{\vect{x}}=0}.
\end{equation*}
To facilitate computations, let us define
$\map{\vect{u}_{\mat{A}}}{\real^n}{\real^m}$ by asking that
\begin{equation*}
u_{\mat{A},a}(\vect{x})=\sign(\inprod{\vect{r}(\mat{A},a)}{\vect{x}}).
\end{equation*}
Note that $B_{\mat{A}}=\vect{u}_{\mat{A}}^{-1}(\vect{0})$ and that on
$\real^n\setminus B_{\mat{A}}$ the function $\vect{u}_{\mat{A}}$ is locally
constant.  Moreover, it is clear that
\begin{equation*}
f(\vect{x})=\inprod{\vect{u}_{\mat{A}}(\vect{x})}
{\mat{A}(\vect{x})}.
\end{equation*}

Now let $\vect{x}_0\in\real^n\setminus B_{\mat{A}}$ be a maximum of $f$
subject to the constraint that $g(\vect{x})=0$\@.  One easily verifies that
$\linder{g}$ has rank $1$ at points that satisfy the constraint.  Thus, by
the Lagrange Multiplier Theorem, there exists $\lambda\in\real$ such that
\begin{equation*}
\linder{(f-\lambda g)}(\vect{x}_0)=\vect{0}.
\end{equation*}
We compute
\begin{equation*}
\linder{f}(\vect{x}_0)\cdot\vect{v}=
\inprod{\vect{u}_{\mat{A}}(\vect{x}_0)}{\mat{A}(\vect{v})},\quad
\linder{g}(\vect{x})\cdot\vect{v}=2\inprod{\vect{x}}{\vect{v}}.
\end{equation*}
Thus $\linder{(f-\lambda g)}(\vect{x}_0)=\vect{0}$ if and only if
\begin{equation*}
\transpose{\mat{A}}(\vect{u}_{\mat{A}}(\vect{x}_0))=2\lambda\vect{x}_0\quad
\implies\quad\snorm{\lambda}=\frac{1}{2}
\dnorm{\transpose{\mat{A}}(\vect{u}_{\mat{A}}(\vect{x}_0))}_2,
\end{equation*}
since $\dnorm{\vect{x}_0}_2=1$\@.  Thus $\lambda=0$ if and only if
$\transpose{\mat{A}}(\vect{u}_{\mat{A}}(\vect{x}_0))=\vect{0}$\@.  Therefore,
if $\lambda=0$\@, then $f(\vect{x}_0)=0$\@.  We can disregard this
possibility since $f$ cannot have a maximum of zero as we are assuming that
$\mat{A}$ has no zero rows.  As $\lambda\not=0$ we have
\begin{equation*}
f(\vect{x}_0)=\inprod{\transpose{\mat{A}}
(\vect{u}_{\mat{A}}(\vect{x}_0))}{\vect{x}_0}=
\frac{1}{2\lambda}\dnorm{\transpose{\mat{A}}
(\vect{u}_{\mat{A}}(\vect{x}_0))}_2^2=2\lambda.
\end{equation*}
We conclude that, at solutions of the constrained maximisation problem, we must
have
\begin{equation*}
f(\vect{x}_0)=\dnorm{\transpose{\mat{A}}(\vect{u})}_2,
\end{equation*}
where $\vect{u}$ varies over the nonzero points in the image of
$\vect{u}_{\mat{A}}$\@,~\ie~over points from $\{-1,1\}^m$\@.

This would conclude the proof of this part of the theorem in the case that
$\mat{A}$ has no zero rows, but for the fact that it is possible that $f$
attains its maximum on $B_{\mat{A}}$\@.  We now show that this does not
happen.  Let $\vect{x}_0\in B_{\mat{A}}$ satisfy $\dnorm{\vect{x}_0}_2=1$ and
denote
\begin{equation*}
A_0=\setdef{a\in\{1,\dots,m\}}{u_{\mat{A},a}(\vect{x}_0)=0}.
\end{equation*}
Let $A_1=\{1,\dots,m\}\setminus A_0$\@.  Let $a_0\in A_0$\@.  For
$\epsilon\in\real$ define
\begin{equation*}
\vect{x}_\epsilon=\frac{\vect{x}_0+\epsilon\vect{r}(\mat{A},a_0)}
{\sqrt{1+\epsilon^2\dnorm{\vect{r}(\mat{A},a_0)}_2^2}}.
\end{equation*}
Note that $\vect{x}_\epsilon$ satisfies the constraint
$\dnorm{\vect{x}_\epsilon}_2^2=1$\@.  Now let $\epsilon_0\in\realp$ be
sufficiently small that
\begin{equation*}
\inprod{\vect{r}(\mat{A},a)}{\vect{x}_\epsilon}\not=0
\end{equation*}
for all $a\in A_1$ and $\epsilon\in\interval[{-\epsilon_0},{\epsilon_0}]$\@.
Then we compute
\begin{align}\notag
\dnorm{\mat{A}(\vect{x}_\epsilon)}_1=&\;
\sum_{a=1}^m\snorm{\inprod{\vect{r}(\mat{A},a)}{\vect{x}_0}+
\epsilon\inprod{\vect{r}(\mat{A},a)}{\vect{r}(\mat{A},a_0)}}+O(\epsilon^2)\\
\notag
=&\;\sum_{a\in A_0}\snorm{\epsilon}
\snorm{\inprod{\vect{r}(\mat{A},a)}{\vect{r}(\mat{A},a_0)}}\\
\label{eq:matnorms1}
&\;+\sum_{a\in A_1}\snorm{\inprod{\vect{r}(\mat{A},a)}
{\vect{x}_0}+\epsilon\inprod{\vect{r}(\mat{A},a)}{\vect{r}(\mat{A},a_0)}}+
O(\epsilon^2).
\end{align}
Since we are assuming that none of the rows of $\mat{A}$ are zero,
\begin{equation}\label{eq:matnorms2}
\sum_{a\in A_0}\snorm{\epsilon}
\snorm{\inprod{\vect{r}(\mat{A},a)}{\vect{r}(\mat{A},a_0)}}>0
\end{equation}
for $\epsilon\in\interval[{-\epsilon_0},{\epsilon_0}]$\@, as long as
$\epsilon_0$ is sufficiently small.  Now take $a\in A_1$\@.  If $\epsilon$ is
sufficiently small we can write
\begin{equation*}
\snorm{\inprod{\vect{r}(\mat{A},a)}
{\vect{x}_0}+\epsilon\inprod{\vect{r}(\mat{A},a)}{\vect{r}(\mat{A},a_0)}}=
\snorm{\inprod{\vect{r}(\mat{A},a)}{\vect{x}_0}}+\epsilon C_a
\end{equation*}
for some $C_a\in\real$\@.  As a result, and using~\eqref{eq:matnorms1}\@, we
have
\begin{equation*}
\dnorm{\mat{A}(\vect{x}_\epsilon)}_1=\dnorm{\mat{A}(\vect{x}_0)}_1+
\sum_{a\in A_0}(\snorm{\epsilon}
\snorm{\inprod{\vect{r}(\mat{A},a)}{\vect{r}(\mat{A},a_0)}}+
\epsilon\sum_{a\in A_1}C_a+O(\epsilon^2).
\end{equation*}
It therefore follows, by choosing $\epsilon_0$ to be sufficiently small, that
we have
\begin{equation*}
\dnorm{\mat{A}(\vect{x}_\epsilon)}_1>\dnorm{\mat{A}(\vect{x}_0)}_1
\end{equation*}
either for all $\epsilon\in\interval[{-\epsilon_0},0)$ or for all
$\epsilon\in\interval(0,{\epsilon_0}]$\@, taking~\eqref{eq:matnorms2} into
account.  Thus if $\vect{x}_0\in B_{\mat{A}}$ then $\vect{x}_0$ is not a
local maximum for $f$ subject to the constraint $g^{-1}(0)$\@.

Finally, suppose that $\mat{A}$ has some rows that are zero.  Let
\begin{equation*}
A_0=\setdef{a\in\{1,\dots,m\}}{\vect{r}(\mat{A},a)=\vect{0}}
\end{equation*}
and let $A_1=\{1,\dots,m\}\setminus A_0$\@.  Let $A_1=\{a_1,\dots,a_k\}$ with
$a_1<\dots<a_k$\@, and define $\hat{\mat{A}}\in\lin{\real^n}{\real^k}$ by
\begin{equation*}
\hat{\mat{A}}(\vect{x})=\sum_{r=1}^k
\inprod{\vect{r}(\mat{A},a_r)}{\vect{x}}\vect{e}_r,
\end{equation*}
and note that $\dnorm{\mat{A}(\vect{x})}_1=\dnorm{\hat{\mat{A}}(\vect{x})}_1$
for every $\vect{x}\in\real^n$\@.  If $\vect{y}\in\real^m$ define
$\hat{\vect{y}}\in\real^k$ by removing from $\vect{y}$ the elements
corresponding to the zero rows of $\mat{A}$\@:
\begin{equation*}
\hat{\vect{y}}=(y_{a_1},\dots,y_{a_k}).
\end{equation*}
Then we easily determine that $\transpose{\mat{A}}(\vect{y})=
\transpose{\hat{\mat{A}}\null}(\hat{\vect{y}})$\@.  Therefore,
\begin{align*}
\dnorm{\mat{A}}_{2,1}=&\;\sup\setdef{\dnorm{\mat{A}(\vect{x})}_1}
{\dnorm{\vect{x}}_2=1}\\
=&\;\sup\setdef{\dnorm{\hat{\mat{A}}(\vect{x})}_1}
{\dnorm{\vect{x}}_2=1}=\dnorm{\hat{\mat{A}}}_{2,1}\\
=&\;\max\setdef{\dnorm{\transpose{\hat{\mat{A}}\null}(\hat{\vect{u}})}_2}
{\hat{\vect{u}}\in\{-1,1\}^k}\\
=&\;\max\setdef{\dnorm{\transpose{\mat{A}}(\vect{u})}_2}
{\vect{u}\in\{-1,1\}^m},
\end{align*}
and this finally gives the result.

\eqref{pl:matnorm22} Note that, in this case, we wish to maximise the function
$\vect{x}\mapsto\dnorm{\mat{A}(\vect{x})}^2_2$ subject to the constraint that
$\dnorm{\vect{x}}^2_2=1$\@.  In this case, the function we are maximising and
the function defining the constraint are infinitely differentiable.
Therefore, we can use the Lagrange Multiplier Theorem to determine the
character of the maxima.  Thus we define
\begin{equation*}
f(\vect{x})=\dnorm{\mat{A}(\vect{x})}_2^2,\quad
g(\vect{x})=\dnorm{\vect{x}}_2^2-1.
\end{equation*}
As $\linder{g}$ has rank $1$ at points satisfying the constraint, if a point
$\vect{x}_0\in\real^n$ solves the constrained maximisation problem, then
there exists $\lambda\in\real$ such that
\begin{equation*}
\linder{(f-\lambda g)}(\vect{x}_0)=0.
\end{equation*}
Since
$f(\vect{x})=\inprod{\transpose{\mat{A}}\scirc\mat{A}(\vect{x})}{\vect{x}}$\@,
we compute
\begin{equation*}
\linder{f}(\vect{x})\cdot\vect{v}=
2\inprod{\transpose{\mat{A}}\scirc\mat{A}(\vect{x})}{\vect{v}}.
\end{equation*}
As above, $\linder{g}(\vect{x})\cdot\vect{v}=2\inprod{\vect{x}}{\vect{v}}$\@.
Thus $\linder{(f-\lambda g)}(\vect{x}_0)=0$ implies that
\begin{equation*}
\transpose{\mat{A}}\scirc\mat{A}(\vect{x}_0)=\lambda\vect{x}_0.
\end{equation*}
Thus it must be the case that $\lambda$ is an eigenvalue for
$\transpose{\mat{A}}\scirc\mat{A}$ with eigenvector $\vect{x}_0$\@.  Since
$\transpose{\mat{A}}\scirc\mat{A}$ is symmetric and positive-semidefinite,
all eigenvalues are real and nonnegative.  Thus there exist
$\lambda_1,\dots,\lambda_n\in\realnn$ and vectors
$\vect{x}_1,\dots,\vect{x}_n$ such that
\begin{equation*}
\lambda_1\le\cdots\le\lambda_n,
\end{equation*}
such that
$\transpose{\mat{A}}\scirc\mat{A}(\vect{x}_j)=\lambda_j\vect{x}_j$\@,
$j\in\{1,\dots,n\}$\@, and such that a solution to the problem of maximising
$f$ with the constraint $g^{-1}(0)$ is obtained by evaluating $f$ at one of
the points $\vect{x}_1,\dots,\vect{x}_n$\@.  Thus the problem can be solved
by evaluating $f$ at this finite collection of points, and determining at
which of these $f$ has its largest value.  A computation gives
$f(\vect{x}_j)=\lambda_j$\@, and this part of the result follows.

\eqref{pl:matnorm2infty} First of all, we note that this part of the theorem
certainly holds when $\mat{A}=\mat{0}$\@.  Thus we shall freely assume that
$\mat{A}$ is nonzero when convenient.  We maximise the function
$\vect{x}\mapsto\dnorm{\mat{A}(\vect{x})}_\infty$ subject to the constraint
that $\dnorm{\vect{x}}^2_2=1$\@.  We shall again use the Lagrange Multiplier
Theorem, defining
\begin{equation*}
f(\vect{x})=\dnorm{\mat{A}(\vect{x})}_\infty,\quad
g(\vect{x})=\dnorm{\vect{x}}_2^2-1.
\end{equation*}
Note that $\mat{A}$ is not differentiable on $\real^n$\@, so we first restrict
to a subset where $f$ is differentiable.  Let us define
\begin{equation*}
\mapdef{S_{\mat{A}}}{\real^n}{\subsets^{\{1,\dots,m\}}}
{\vect{x}}{\setdef{a\in\{1,\dots,m\}}
{\inprod{\vect{r}(\mat{A},a)}{\vect{x}}=\dnorm{\mat{A}(\vect{x})}_\infty}.}
\end{equation*}
Then denote
\begin{equation*}
B_{\mat{A}}=\setdef{\vect{x}\in\real^n}{\card(S_{\mat{A}}(\vect{x}))>1}.
\end{equation*}
We easily see that $f$ is differentiable at points that are not in the set
$B_{\mat{A}}$\@.

Let us first suppose that $\vect{x}_0\in\real^n\setminus B_{\mat{A}}$ is a
maximum of $f$ subject to the constraint that $g(\vect{x})=0$\@.  Then there
exists a unique $a_0\in\{1,\dots,m\}$ such that
$f(\vect{x}_0)=\inprod{\vect{r}(\mat{A},a_0)}{\vect{x}_0}$\@.  Since we are
assuming that $\mat{A}$ is nonzero, it must be that $\vect{r}(\mat{A},a_0)$
is nonzero.  Moreover, there exists a neighbourhood $U$ of $\vect{x}_0$ such
that
\begin{equation*}
\sign(\inprod{\vect{r}(\mat{A},a_0)}{\vect{x}})=
\sign(\inprod{\vect{r}(\mat{A},a_0)}{\vect{x}_0})
\end{equation*}
and $f(\vect{x})=\inprod{\vect{r}(\mat{A},a_0)}{\vect{x}}$ for each
$\vect{x}\in U$\@.  Abbreviating
\begin{equation*}
u_{\mat{A},a_0}(\vect{x})=
\sign(\inprod{\vect{r}(\mat{A},a_0)}{\vect{x}}),
\end{equation*}
we have
\begin{equation*}
f(\vect{x})=u_{\mat{A},j}(\vect{x}_0)\inprod{\vect{r}(\mat{A},a_0)}{\vect{x}}
\end{equation*}
for every $\vect{x}\in U$\@.  Note that, as in the proofs of
parts~\eqref{pl:matnorm21} and~\eqref{pl:matnorm22} above,
$\linder{g}(\vect{x})$ has rank $1$ for $\vect{x}\not=0$\@.  Therefore there
exists $\lambda\in\real$ such that
\begin{equation*}
\linder{(f-\lambda g)}(\vect{x}_0)=\vect{0}.
\end{equation*}
We compute
\begin{equation*}
\linder{(f-\lambda g)}(\vect{x}_0)\cdot\vect{v}=
u_{\mat{A},j}(\vect{x}_0)\inprod{\vect{r}(\mat{A},a_0)}{\vect{v}}-
2\lambda\inprod{\vect{x}_0}{\vect{v}}
\end{equation*}
for every $\vect{v}\in\real^n$\@.  Thus we must have
\begin{equation*}
2\lambda\vect{x}_0=u_{\mat{A},a_0}(\vect{x}_0)\vect{r}(\mat{A},a_0).
\end{equation*}
This implies that $\vect{x}_0$ and $\vect{r}(\mat{A},a_0)$ are collinear and
that
\begin{equation*}
\snorm{\lambda}=\frac{1}{2}\dnorm{\vect{r}(\mat{A},a_0)}_2
\end{equation*}
since $\dnorm{\vect{x}_0}_2=1$\@.  Therefore,
\begin{equation*}
f(\vect{x}_0)=u_{\mat{A},a_0}(\vect{x}_0)\inprod{\vect{r}(\mat{A},a_0)}
{\tfrac{1}{2\lambda}u_{\mat{A},a_0}(\vect{x}_0)\vect{r}(\mat{A},a_0)}=
2\lambda.
\end{equation*}
Since $\snorm{\lambda}=\frac{1}{2}\dnorm{\vect{r}(\mat{A},a_0)}_2$ it follows
that
\begin{equation*}
f(\vect{x}_0)=\dnorm{\vect{r}(\mat{A},a_0)}_2.
\end{equation*}

This completes the proof, but for the fact that maxima of $f$ may occur at
points in $B_{\mat{A}}$\@.  Thus let $\vect{x}_0\in B_{\mat{A}}$ be such that
$\dnorm{\vect{x}_0}_2=1$\@.  For $a\in S_{\mat{A}}(\vect{x}_0)$ let us write
\begin{equation*}
\vect{r}(\mat{A},a)=\rho_a\vect{x}_0+\vect{y}_a,
\end{equation*}
where $\inprod{\vect{x}_0}{\vect{y}_a}=0$\@.  Therefore,
$\inprod{\vect{r}(\mat{A},a)}{\vect{x}_0}=\rho_a$\@.  We claim that if there
exists $a_0\in S_{\mat{A}}(\vect{x}_0)$ for which
$\vect{y}_{a_0}\not=\vect{0}$\@, then $\vect{x}_0$ cannot be a maximum of $f$
subject to the constraint $g^{-1}(0)$\@.  Indeed, if $\vect{y}_{a_0}\not=0$
then define
\begin{equation*}
\vect{x}_\epsilon=\frac{\vect{x}_0+\epsilon\vect{y}_{a_0}}
{\sqrt{1+\epsilon^2\dnorm{\vect{y}_{a_0}}_2^2}}.
\end{equation*}
As in the proof of part~\eqref{pl:matnorm21} above, one shows that
$\vect{x}_\epsilon$ satisfies the constraint for every $\epsilon\in\real$\@.
Also as in the proof of part~\eqref{pl:matnorm21}\@, we have
\begin{equation*}
\vect{x}_\epsilon=\vect{x}_0+\epsilon\vect{y}_0+O(\epsilon^2).
\end{equation*}
Thus, for $\epsilon$ sufficiently small,
\begin{equation*}
\snorm{\inprod{\vect{r}(\mat{A},a_0)}{\vect{x}_\epsilon}}=
\snorm{\inprod{\vect{r}(\mat{A},a_0)}{\vect{x}_0}}+
\epsilon C_{a_0}+O(\epsilon^2)
\end{equation*}
where $C_{a_0}$ is nonzero.  Therefore, there exists $\epsilon_0\in\realp$
such that
\begin{equation*}
\snorm{\inprod{\vect{r}(\mat{A},a_0)}{\vect{x}_\epsilon}}>
\snorm{\inprod{\vect{r}(\mat{A},a_0)}{\vect{x}_0}}
\end{equation*}
either for all $\epsilon\in\interval[{-\epsilon_0},0)$ or for all
$\epsilon\in\interval(0,{\epsilon_0}]$\@.  In either case, $\vect{x}_0$
cannot be a maximum for $f$ subject to the constraint $g^{-1}(0)$\@.

Finally, suppose that $\vect{x}_0\in B_{\mat{A}}$ is a maximum for $f$ subject
to the constraint $g^{-1}(0)$\@.  Then, as we saw in the preceding paragraph,
for each $a\in S_{\mat{A}}(\vect{x}_0)$\@, we must have
\begin{equation*}
\vect{r}(\mat{A},a)=\inprod{\vect{r}(\mat{A},a)}{\vect{x}_0}\vect{x}_0.
\end{equation*}
It follows that $\dnorm{\vect{r}(\mat{A},a)}_2^2=
\inprod{\vect{r}(\mat{A},a)}{\vect{x}_0}^2$\@.  Moreover, by definition of
$S_{\mat{A}}(\vect{x}_0)$ and since we are supposing that $\vect{x}_0$ is a
maximum for $f$ subject to the constraint $g^{-1}(0)$\@, we have
\begin{equation}\label{eq:matnorms3}
\dnorm{\vect{r}(\mat{A},a)}_2=\dnorm{\mat{A}}_{2,\infty}.
\end{equation}
Now, if $a\in\{1,\dots,m\}$\@, we claim that
\begin{equation}\label{eq:matnorms4}
\dnorm{\vect{r}(\mat{A},a)}_2\le\dnorm{\mat{A}}_{2,\infty}.
\end{equation}
Indeed suppose that $a\in\{1,\dots,m\}$ satisfies
\begin{equation*}
\dnorm{\vect{r}(\mat{A},a)}_2>\dnorm{\mat{A}}_{2,\infty}.
\end{equation*}
Define $\vect{x}=\frac{\vect{r}(\mat{A},a)}{\dnorm{\vect{r}(\mat{A},a)}_2}$ so
that $\vect{x}$ satisfies the constraint $g(\vect{x})=0$\@.  Moreover,
\begin{equation*}
f(\vect{x})\ge\inprod{\vect{r}(\mat{A},a)}{\vect{x}}=
\dnorm{\vect{r}(\mat{A},a)}_2>\dnorm{\mat{A}}_{2,\infty},
\end{equation*}
contradicting the assumption that $\vect{x}_0$ is a maximum for $f$\@.  Thus,
given that~\eqref{eq:matnorms3} holds for every $a\in S_{\mat{A}}(\vect{x}_0)$
and~\eqref{eq:matnorms4} holds for every $a\in\{1,\dots,m\}$\@, we have
\begin{equation*}
\dnorm{\mat{A}}_{2,\infty}=\max\setdef{\dnorm{\vect{r}(\mat{A},a)}_2}
{a\in\{1,\dots,m\}},
\end{equation*}
as desired.

For the last three parts of the theorem, the following result is useful.
\begin{procnonum}[Lemma]
Let\/ $\dnorm{\cdot}$ be a norm on\/ $\real^n$ and let\/
$\tnorm{\cdot}_\infty$ be the norm induced on\/ $\lin{\real^n}{\real^m}$ by
the norm\/ $\dnorm{\cdot}_\infty$ on\/ $\real^n$ and the norm\/
$\dnorm{\cdot}$ on\/ $\real^m$\@.  Then
\begin{equation*}
\tnorm{\mat{A}}_\infty=\max\setdef{\dnorm{\mat{A}(\vect{u})}}
{\vect{u}\in\{-1,1\}^n}.
\end{equation*}
\begin{subproof}
Note that the set
\begin{equation*}
\setdef{\vect{x}\in\real^n}{\dnorm{\vect{x}}_\infty\le1}
\end{equation*}
is a convex polytope.  Therefore, this set is the convex hull of the vertices
$\{-1,1\}^n$\@; see~\cite[Theorem~2.6.16]{RJW:94}\@.  Thus, if
$\dnorm{\vect{x}}_\infty=1$ we can write
\begin{equation*}
\vect{x}=\sum_{\vect{u}\in\{-1,1\}^n}\lambda_{\vect{u}}\vect{u}
\end{equation*}
where $\lambda_{\vect{u}}\in\interval[0,1]$ for each $\vect{u}\in\{-1,1\}^n$
and
\begin{equation*}
\sum_{\vect{u}\in\{-1,1\}^n}\lambda_{\vect{u}}=1.
\end{equation*}
Therefore,
\begin{align*}
\dnorm{\mat{A}(\vect{x})}=&\;
\adnorm{\sum_{\vect{u}\in\{-1,1\}^n}\lambda_{\vect{u}}\mat{A}(\vect{u})}\le
\sum_{\vect{u}\in\{-1,1\}^n}\lambda_{\vect{u}}
\dnorm{\mat{A}(\vect{u})}\\
\le&\;\left(\sum_{\vect{u}\in\{-1,1\}^n}\lambda_{\vect{u}}\right)
\max\setdef{\dnorm{\mat{A}(\vect{u})}}{\vect{u}\in\{-1,1\}^n}\\
=&\;\max\setdef{\dnorm{\mat{A}(\vect{u})}}{\vect{u}\in\{-1,1\}^n}.
\end{align*}
Therefore,
\begin{equation*}
\sup\setdef{\dnorm{\mat{A}(\vect{x})}}{\dnorm{\vect{x}}_\infty=1}\le
\max\setdef{\dnorm{\mat{A}(\vect{u})}}{\vect{u}\in\{-1,1\}^n}
\le\sup\setdef{\dnorm{\mat{A}(\vect{x})}}{\dnorm{\vect{x}}_\infty=1},
\end{equation*}
the last inequality holding since if $\vect{u}\in\{-1,1\}^n$ then
$\dnorm{\vect{u}}_\infty=1$\@.  The result follows since the previous
inequalities must be equalities.
\end{subproof}
\end{procnonum}

\eqref{pl:matnorminfty1} This follows immediately from the preceding lemma.

\eqref{pl:matnorminfty2} This too follows immediately from the preceding
lemma.

\eqref{pl:matnorminftyinfty} Note that for $\vect{u}\in\{-1,1\}^n$ we have
\begin{equation*}
\snorm{\inprod{\vect{r}(\mat{A},a)}{\vect{u}}}=
\asnorm{\sum_{j=1}^nA_{aj}u_j}\le\sum_{j=1}^n\snorm{A_{aj}}=
\dnorm{\vect{r}(\mat{A},a)}_1.
\end{equation*}
Therefore, using the previous lemma,
\begin{align*}
\dnorm{\mat{A}}_{\infty,\infty}=&\;\max\setdef{\dnorm{\mat{A}(\vect{u})}_\infty}
{\vect{u}\in\{-1,1\}^n}\\
=&\;\max\setdef{\max\setdef{\snorm{\inprod{\vect{r}(\mat{A},a)}{\vect{u}}}}
{a\in\{1,\dots,m\}}}{\vect{u}\in\{-1,1\}^n}\\
\le&\;\max\setdef{\dnorm{\vect{r}(\mat{A},a)}_1}{a\in\{1,\dots,m\}}.
\end{align*}
To establish the other inequality, for $a\in\{1,\dots,m\}$ define
$\vect{u}_a\in\{-1,1\}^n$ by
\begin{equation*}
u_{a,j}=\begin{cases}1,&A_{aj}\ge0,\\-1,&A_{aj}<0\end{cases}
\end{equation*}
and note that a direct computation gives the $a$th component of
$\mat{A}(\vect{u}_a)$ as $\dnorm{\vect{r}(\mat{A},a)}_1$\@.  Therefore,
\begin{align*}
\max\setdef{\dnorm{\vect{r}(\mat{A},a)}_1}{a\in\{1,\dots,m\}}=&\;
\max\setdef{\snorm{\mat{A}(\vect{u}_a)_a}}{a\in\{1,\dots,m\}}\\
\le&\;\max\setdef{\dnorm{\mat{A}(\vect{u})}_\infty}{\vect{u}\in\{-1,1\}^n}=
\dnorm{\mat{A}}_{\infty,\infty},
\end{align*}
giving this part of the theorem.
\end{proof}
\end{theorem}

\section{Complexity of induced norm computations}

Let us consider a comparison of the nine induced matrix norms in terms of the
computational effort required.  One would like to know how many operations
are required to compute any of the norms.  We shall do this making the
following assumptions on our computational model.
\begin{quote}
Floating point operations are carried out to an accuracy of
$\epsilon=2^{-N}$ for some fixed $N\in\integerp$\@.  By $M(N)$ we denote the
number of operations required to multiply integers $j_1$ and $j_2$ satisfying
$0\le j_1,j_2\le 2^N$\@.  We assume that addition and multiplication of
floating point numbers can be performed with a relative error of $O(2^{-N})$
using $O(M(N))$ operations.
\end{quote}
With this assumption, we can deduce the computational complexity of the basic
operations we will need.
\begin{compactenum}
\item Computing a square root takes $O(M(N))$ operations;
see~\cite{RPB:76}\@.
\item Computing the absolute value of a number is $1$ operation (a bit flip).
\item Comparing two numbers takes $O(N)$ operations.
\item Finding the maximum number in a list of $k$ numbers takes $O(kN)$
operations; see~\cite{MB/RWF/VP/RLR/RET:73}\@.
\item If $\mat{A}\in\lin{\real^n}{\real^m}$ and
$\mat{B}\in\lin{\real^m}{\real^p}$ then the matrix multiplication
$\mat{B}\mat{A}$ takes $O(mnpM(n))$ operations.  Faster matrix multiplication
algorithms are possible than the direct one whose complexity we describe
here,~\cite[\eg][]{DC/SW:90}\@, but we are mainly interested in the fact that
matrix multiplication has polynomial complexity in the size of the matrices.
\item Computation of the $QR$-decomposition of an $k\times k$ matrix
$\mat{A}$ has computational complexity $O(k^3M(N))$\@.  Note that the
$QR$-decomposition can be used to determine the Gram\textendash{}Schmidt
orthogonalisation of a finite number of vectors.  We refer
to~\cite[\S5.2]{GHG/CFVL:96} for details.
\item Let us describe deterministic bounds for the operations needed to
compute the eigenvalues and eigenvectors of a $k\times k$ matrix $\mat{A}$\@,
following \citet{VYP/ZQC:99}\@.  Let us fix some norm $\tnorm{\cdot}$ on
$\lin{\real^k}{\real^k}$\@.  Given $\epsilon=2^{-N}$ as above, let
$\beta\in\realp$ be such that $2^{-\beta}\tnorm{\mat{A}}\le\epsilon$\@.  Then
\citeauthor{VYP/ZQC:99} show that the eigenvalues and eigenvectors of
$\mat{A}$ can be computed for $\mat{A}$ using an algorithm of complexity
\begin{equation}\label{eq:eigencomplexity}
O(k^3M(N))+O((k\log^2k)(\log\beta+\log^2k)M(N)).
\end{equation}
There are stochastic, iterative, or gradient flow algorithms that will
generically perform computations with fewer operations than predicted by this
bound.  However, the complexity of such algorithms is difficult to
understand, or they require unbounded numbers of operations in the worst
case.  In any event, here we only care that the complexity of the
eigenproblem is polynomial.
\item The previous two computational complexity results can be combined to
show that finding the square root of a symmetric positive-definite matrix has
computational complexity given by~\eqref{eq:eigencomplexity}\@.  This is no
doubt known, but let us see how this works since it is simple.  First compute
the eigenvalues and eigenvectors of $\mat{A}$ using an algorithm with
complexity given by~\eqref{eq:eigencomplexity}\@.  The eigenvectors can be
made into an orthonormal basis of eigenvectors using the
Gram\textendash{}Schmidt procedure.  This decomposition can be performed
using an algorithm of complexity $O(n^3M(N))$\@.  Assembling the orthogonal
eigenvectors into the columns of a matrix gives an orthogonal matrix
$\mat{U}\in\lin{\real^n}{\real^n}$ and a diagonal matrix
$\mat{D}\in\lin{\real^n}{\real^n}$ with positive diagonals such that
$\mat{A}=\mat{U}\mat{D}\transpose{\mat{U}}$\@.  Then the
matrix $\mat{D}^{1/2}$ with diagonal entries equal to the square roots of the
diagonal of $\mat{D}$ can be constructed with complexity $O(nM(n))$\@.
Finally, $\mat{A}^{1/2}=\mat{U}\mat{D}^{1/2}\mat{U}$ is computed using matrix
multiplication with complexity of $O(n^3M(n))$\@.
\end{compactenum}

Using these known computational complexity results, it is relatively
straightforward to assess the complexity of the computations of the various
norms in Theorem~\ref{the:matnorms}\@.  In Table~\ref{tab:complexity}
\begin{table}[htbp]
\centering
\caption{Complexity of computing the norms
$\dnorm{\cdot}_{p,q}$}\label{tab:complexity}
\begin{tabular}{c|ccc}
\toprule
\diagbox{$p$}{$q$}&$1$&$2$&$\infty$\\\midrule
$1$&$O(mn)$&$O(mn)$&$O(mn)$\\
$2$&$O(mn2^m)$&$O(n^3)$&$O(mn)$\\
$\infty$&$O(mn2^n)$&$O(mn2^n)$&$O(mn)$\\\bottomrule
\end{tabular}
\end{table}
we display this data, recording only the dependency of the computations on
the number of rows $m$ and columns $n$ of the matrix.  Note that the cases of
$(p,q)\in\{(2,1),(\infty,1),(\infty,2)\}$ are exceptional in that the
required operations grow exponentially with the size of $\mat{A}$\@.  One
must exercise some care in drawing conclusions here.  For example, as we show
in the proof of Theorem~\ref{the:matnorms}\@,
\begin{equation}\label{eq:inftyinfty}
\dnorm{\mat{A}}_{\infty,\infty}=\max\setdef{\dnorm{\mat{A}(\vect{u})}_\infty}
{\vect{u}\in\{-1,1\}^n},
\end{equation}
and this computation has complexity $O(mn2^n)$\@.  However, it turns out that
the norm can be determined with a formula that is actually less complex.
Indeed, our proof of the formula for
$\dnorm{\cdot}_{\infty,\infty}$\@\textemdash{}which is not the usual
proof\textemdash{}starts with the formula~\eqref{eq:inftyinfty} and produces a
result with complexity $O(mn)$ as stated in Table~\ref{tab:complexity}\@.

One is then led to ask, are there similar simplifications of the norms
corresponding to the cases $(p,q)\in\{(2,1),(\infty,1),(\infty,2)\}$\@?
\citet{JR:00} shows that the computation of $\dnorm{\cdot}_{\infty,1}$ is
NP-hard.  We shall show here, using his ideas, that the computation of the
norms $\dnorm{\cdot}_{2,1}$ and $\dnorm{\cdot}_{\infty,2}$ are likewise
difficult, perhaps impossible, to reduce to algorithms with polynomial
complexity.
\begin{theorem}
If there exists an algorithm to compute\/ $\dnorm{\mat{A}}_{2,1}$ or\/
$\dnorm{\mat{A}}_{\infty,2}$ whose computational complexity is polynomial in
the number of rows and the number of columns of\/ $\mat{A}$\@, then P=NP.
\begin{proof}
First note that
$\dnorm{\mat{A}}_{2,1}=\dnorm{\transpose{\mat{A}}}_{\infty,2}$\@, so it
suffices to prove the theorem only for $(p,q)=(\infty,2)$\@.

Following \citet{JR:00} we introduce the notion of an \defn{$MC$-matrix}
(``$MC$'' stands for ``max-cut'' since these matrices are related to the
``max-cut problem'' in graph theory) as a symmetric matrix
$\mat{A}\in\lin{\real^n}{\real^n}$ with the property that the diagonal
elements are equal to $n$ and the off-diagonal elements are either $0$ or
$-1$\@.  \cite{JR:94} shows that $MC$-matrices are positive-definite.
\citet{SP/JR:93} also prove the following.
\begin{quote}\em
The following decision problem is NP-complete:\\
Given an\/ $n\times n$\/ $MC$-matrix\/ $\mat{A}$ and\/ $M\in\integerp$\@,
is\/ $\inprod{\mat{A}(\vect{u})}{\vect{u}}\ge M$ for some\/
$\vect{u}\in\{-1,1\}^n$?
\end{quote}
We will use this fact crucially in our proof.

Let us call a symmetric matrix $\mat{A}\in\lin{\real^n}{\real^n}$ a
\defn{$\sqrt{MC}$-matrix} if $\mat{A}\scirc\mat{A}$ is an $MC$-matrix.  Note
that the map $\mat{A}\mapsto\mat{A}\scirc\mat{A}$ from the set of
$\sqrt{MC}$-matrices to the set of $MC$-matrices is surjective since
$MC$-matrices have symmetric positive-definite square roots by virtue of
their being themselves symmetric and positive-definite.

Now suppose that there exists an algorithm for determining the
$(\infty,2)$-norm of a matrix, the computational complexity of which is of
polynomial order in the number of rows and columns of the matrix.  Let
$\mat{A}$ be an $n\times n$ $MC$-matrix and let $M\in\integerp$\@.  As we
pointed out prior to stating the theorem, one can determine the
$\sqrt{MC}$-matrix $\mat{A}^{1/2}$ using an algorithm with computational
complexity that is polynomial in $n$\@.  Then, by assumption, we can compute
\begin{align*}
\dnorm{\mat{A}^{1/2}}_{\infty,2}^2=&\;
\max\setdef{\dnorm{\mat{A}^{1/2}(\vect{u})}_2^2}{\vect{u}\in\{-1,1\}^n}\\
=&\;\max\setdef{\inprod{\mat{A}(\vect{u})}{\vect{u}}}{\vect{u}\in\{-1,1\}^n}
\end{align*}
in polynomial time.  In particular, we can determine whether
$\inprod{\mat{A}(\vect{u})}{\vect{u}}\ge M$ in polynomial time.  As we stated
above, this latter decision problem is NP-complete, and so we must have P=NP.
\end{proof}
\end{theorem}

\printbibliography[heading=bibintoc]

\end{document}